%% file: main.tex
\pdfoutput=1
\newif\ifSpringer

\Springerfalse 

\ifSpringer

\RequirePackage{fix-cm}
\documentclass[smallextended]{svproc}
\usepackage{etoolbox}
\usepackage[misc]{ifsym}
\AtEndEnvironment{proof}{\phantom{}\qed}
\else
\documentclass[11pt]{scrartcl}
\usepackage{amsthm}

\fi

\usepackage[utf8]{inputenc}
\usepackage[T1]{fontenc}
\usepackage{lmodern}	
\usepackage[english]{babel}

\usepackage[draft=false,kerning=true]{microtype}
\usepackage[numbers]{natbib} 
\usepackage{url}

\usepackage{amsmath,amssymb,graphicx}
\usepackage{pdflscape}
\usepackage{fancyhdr}
\usepackage{adjustbox}
\usepackage{rotating}
\usepackage{float}
\usepackage[caption = false]{subfig}
\usepackage[ruled,vlined,linesnumbered]{algorithm2e}
\usepackage{booktabs,tabularx} 
\usepackage{siunitx} 
\usepackage{pgfplots}
\usepackage{pgfplotstable} 
\usepackage{longtable,cellspace}
\usepackage{subfig}
\usepackage{tablefootnote}

\usepackage[title]{appendix}

\usepackage[colorinlistoftodos,prependcaption]{todonotes}

\usepackage{scalerel}
\usepackage{tikz}
\usetikzlibrary{svg.path}

\addparagraphcolumntypes{X}

\DeclareOldFontCommand{\bf}{\normalfont\bfseries}{\mathbf}

\newcommand{\R}{\mathbb{R}}
\newcommand{\sn}{{\mathcal S}^n}
\newcommand{\inprod}[2]{{\langle #1,#2 \rangle}} 
\newcommand{\trace}{\textrm{trace}}

\newcommand{\nosemic}{\renewcommand{\@endalgocfline}{\relax}}
\newcommand{\dosemic}{\renewcommand{\@endalgocfline}{\algocf@endline}}
\let\oldnl\nl
\newcommand{\nonl}{\renewcommand{\nl}{\let\nl\oldnl}}

\definecolor{orcidlogocol}{HTML}{A6CE39}
\tikzset{
	orcidlogo/.pic={
		\fill[orcidlogocol] svg{M256,128c0,70.7-57.3,128-128,128C57.3,256,0,198.7,0,128C0,57.3,57.3,0,128,0C198.7,0,256,57.3,256,128z};
		\fill[white] svg{M86.3,186.2H70.9V79.1h15.4v48.4V186.2z}
		svg{M108.9,79.1h41.6c39.6,0,57,28.3,57,53.6c0,27.5-21.5,53.6-56.8,53.6h-41.8V79.1z M124.3,172.4h24.5c34.9,0,42.9-26.5,42.9-39.7c0-21.5-13.7-39.7-43.7-39.7h-23.7V172.4z}
		svg{M88.7,56.8c0,5.5-4.5,10.1-10.1,10.1c-5.6,0-10.1-4.6-10.1-10.1c0-5.6,4.5-10.1,10.1-10.1C84.2,46.7,88.7,51.3,88.7,56.8z};
	}
}
\providecommand{\keywords}[1]
{
	\small	
	\textbf{\textit{Keywords---}} #1
}

\newcommand\orcidicon[1]{\href{https://orcid.org/#1}{\mbox{\scalerel*{
				\begin{tikzpicture}[yscale=-1,transform shape]
				\pic{orcidlogo};
				\end{tikzpicture}
			}{|}}}}
\usepackage{hyperref}

\fancypagestyle{mylandscape}{
	\fancyhf{} 
	\fancyfoot{
		\makebox[\textwidth][r]{
			\rlap{\hspace{.20cm}
				\smash{
					\raisebox{4.87in}{
						\rotatebox{90}{\thepage}}}}}}
}

\pgfplotstableset{
	begin table=\begin{longtable},
		end table=\end{longtable},
}

\pgfplotstableset{
	highlight col max/.code 2 args={
		\findmax{#1}{#2}{\maxval}
		\edef\setstyles{\noexpand\pgfplotstableset{
				every row \maxval\noexpand\space column #2/.style={
					postproc cell content/.append style={
						/pgfplots/table/@cell content/.add={$\noexpand\bf}{$}
					},
				}
			}
		}\setstyles
	},
	highlight col min/.code 2 args={
		\findmin{#1}{#2}{\minval}
		\edef\setstyles{\noexpand\pgfplotstableset{
				every row \minval\noexpand\space column #2/.style={
					postproc cell content/.append style={
						/pgfplots/table/@cell content/.add={\noexpand\color{red}$\noexpand\bf}{$}
					},
				}
			}
		}\setstyles
	},
	highlight row max/.code 2 args={
		\pgfmathtruncatemacro\rowindex{#2-1}
		\pgfplotstabletranspose{\transposed}{#1}
		\findmax{\transposed}{\rowindex}{\maxval}
		\edef\setstyles{\noexpand\pgfplotstableset{
				every row \rowindex\space column \maxval\noexpand/.style={
					postproc cell content/.append style={
						/pgfplots/table/@cell content/.add={$\noexpand\bf}{$}
					},
				}
			}
		}\setstyles
	},
	highlight row min/.code 2 args={
		\pgfmathtruncatemacro\rowindex{#2-1}
		\pgfplotstabletranspose{\transposed}{#1}
		\findmin{\transposed}{\rowindex}{\maxval}
		\edef\setstyles{\noexpand\pgfplotstableset{
				every row \rowindex\space column \maxval\noexpand/.style={
					skip cols between index={0}{1},
					postproc cell content/.append style={
						/pgfplots/table/@cell content/.add={\noexpand\color{red}$\noexpand\bf}{$}
					},
				}
			}
		}\setstyles
	},
}
\makeatletter
\long\def\pgfplotstabletypeset@opt@collectarg[#1]#2{%
	
	\pgfplotstable@isloadedtable{#2}%
	{\pgfplotstabletypeset@opt@[#1]{#2}}%
	{\pgfplotstabletypesetfile@opt@[#1]{#2}}%
}
\makeatother

\begin{document}

\title{Tight SDP relaxations for cardinality-constrained problems \thanks{This project has received funding from the European Union’s Horizon 2020 research and innovation programme under the Marie Sk\l{}odowska-Curie grant agreement MINOA No 764759.}
}

%

\ifSpringer
\else

\author{Angelika Wiegele \orcidicon{0000-0003-1670-7951} 
	\footnote{Institut f\"ur Mathematik, Alpen-Adria-Universit\"at Klagenfurt, Universit\"atsstraße 65-67, 9020, \href{mailto:angelika.wiegele@aau.at}{angelika.wiegele@aau.at}, \href{mailto:shudian.zhao@aau.at}{shudian.zhao@aau.at}} \,
	\and Shudian Zhao \orcidicon{0000-0001-6352-0968} $^\dagger$\footnote{Corresponding author}
}

\fi


\maketitle

\begin{abstract}
We model the cardinality-constrained portfolio problem using semidefinite matrices and investigate a relaxation using semidefinite programming. Experimental results show that this relaxation generates tight lower bounds and even achieves optimality on many instances from the literature. This underlines the modeling power of semidefinite programming for mixed-integer quadratic problems. 

\keywords{Semidefinite programming, cardinality-constrained problem, mixed-integer nonlinear programming}
\end{abstract}

\section{Introduction}

The cardinality-constrained optimization problem is widely applied in areas of finance such as the portfolio optimization problem where the number of stocks to be traded is bounded. 
It is even NP-complete to test the feasibility of an optimization problem with cardinality constraints \cite{bienstock1996computational}.

There are various previous studies about cardinality-constrained problems. \citet{frangioni2007sdp} and \citet{zheng2014improving} work with the approaches that use diagonalizations and perspective cuts to improve the  continuous relaxation bound in order to deal with the cardinality-constrained problems with lower bounds for nonzero stocks. \citet{burdakov2016mathematical} have  proposed nonconvex relaxations that still have the same solutions (in the sense of global minima), and then tools for minimization problems in continuous variables are applied to solve the relaxations. This problem also lies in the family of quadratic programming problems with complementary constraints, where \citet{braun2005semidefinite} have introduced heuristics via semidefinite programming to generate upper bounds for this problem.

Semidefinite programming belongs to the field of convex optimization. A semidefinite programming (SDP) problem can be solved in polynomial time. SDP has shown advantages in approximating integer or mixed-integer nonlinear programming problems. In particular, quadratic problems can be quite naturally modeled by semidefinite programming. For quadratically constrained quadratic programming (QCQP), the optimal solution of the SDP relaxation is tight when the QCQP problem is convex~\cite{park2017general}. Interior point methods (IPMs) are commonly used to solve SDP problems and implemented in SDP solvers such as Mosek \cite{mosek}. IPMs can solve SDP problems to high precision as long as the number of constraints and the size of the matrices is reasonable.

\subsubsection*{Notation} The notation $[n]$ stands for the set of integers $\{1,\dots,n\}$.
We denote by $\sn$ the set of symmetric $n\times n$ matrices. The operation $\textrm{diag}(M)$ maps the diagonal entries of matrix $M$ to a vector. We denote by $\inprod{\cdot}{\cdot}$ the trace inner product. That is, for any $M, N \in \R^{m\times n}$, we define $\inprod{M}{N}:= \trace (N^\top M)$.
We write $X\succeq Y$ if matrix $X-Y$ is positive semidefinite. 


\section{The cardinality-constrained portfolio optimization problem}
Given $n$ stocks and the covariance matrix $Q$, a profit vector $\mu$ and the minimum expected return $\rho$. The objective is to find a portfolio that minimizes the risk while a minimum expected return is achieved. A mathematical programming formulation is as follows.
\begin{equation}
	\begin{aligned}
	\min~&x^\top Q x\\
		\textrm{s.t.}~&\mu^\top x \geq \rho,\\
		&e^\top x = 1,\\
		& 0 \leq x_i\leq u_i,~\forall i \in [n],
	\end{aligned}
\end{equation}
where $x$ indicates the weights for each stock of a portfolio which is nonnegative and bounded by $u$, $x^\top Q x $ is the standard deviation of the portfolio, which is used to measure the risk of that portfolio, and $\mu^\top x $ is the expected return.

One commonly used constraint for a portfolio problem is the cardinality constraint, in other words, the maximum number of stocks to be chosen.
\begin{equation}\label{eq:card}
\begin{aligned}
		\min~&x^\top Q x\\
	\textrm{s.t.}~&\mu^\top x \geq \rho,\\
	& e^\top x = 1,\\
	& 0 \leq x_i\leq u_i,~\forall i \in [n],\\
	& \textbf{Card}(x) \leq \aleph,
	\end{aligned}
\end{equation}
where $\textbf{Card}(x)$ indicates the number of nonzero components of $x$.

Problem \eqref{eq:card} can be formed as a mixed-integer nonlinear programming problem as follows~\cite{burdakov2016mathematical}.
\begin{equation}\label{eq:cardMIQP}
\begin{aligned}
    \min_{x,y}~&x^\top Qx \\ 
    \textrm{s.t.}~&\mu^\top x \geq \rho,\\
    &e^\top x \leq 1,\\
    &e^\top y \geq n -\aleph,\\
    & x_i y_i = 0,~\forall i \in [n],\\
    & y_i \in \{ 0, 1\},~\forall i \in [n],\\
    & 0 \leq x_i \leq u_i, ~\forall i \in [n],
\end{aligned}
\end{equation}
where the slack variables $y\in \{0,1\} ^n$ are integer and the complementary constraints $x_i y_i =0$ relates $x$ and $y$. Hence, $x_i$ and $y_i$ cannot be both positive, and the constraint $e^\top y \geq n - \aleph$ requires that $y$ has at least $n-\aleph$ nonzero elements, in return, then $x$ has at most $\aleph$ nonzero elements.

This problem is NP-hard \cite{judice1994linear}. When $\aleph$ is small, an optimal solution can be found easily by global search, while it gets more complicated to solve when $\aleph$ increases.

\section{A Semidefinite Programming Relaxation}

Before introducing a relaxation using semidefinite programming, note that in problem~\eqref{eq:cardMIQP} we have $x^\top Qx = \langle x, Qx \rangle = \langle Q, xx^\top \rangle$. Moreover, $y_i \in \{0,1\}$ implies $y_i^2 = y_i$ and hence the diagonal of matrix $yy^\top$ must be equal to $y$.
The condition $x_iy_i=0$ translates to the diagonal of matrix $xy^\top$ being $0$. 

We now introduce matrices $X$ and $Y$, substitute the term $xx^\top$ by the symmetric matrix $X$, i.e., $X = xx^\top$,
and relax this condition to 
$X \succeq x x^\top$. Similarly, we relax $Y = yy^\top$ to $Y \succeq y y^\top$.

Then, with the Schur complement we have 
\begin{equation*}
\begin{aligned}
X \succeq x x^\top,~ Y \succeq y y^\top \iff \exists  Z \in \R^{n\times n} ~ \textrm{s.t.}~\begin{pmatrix} 1 & x^\top & y^\top \\ x & X & Z \\ y & Z^\top & Y  \end{pmatrix} \succeq 0.
\end{aligned}
\end{equation*}
Hence, the final SDP relaxation is given as
\begin{equation} \label{eq:SDP} 
 \begin{aligned}
        \min_{\bar{X}}~ & \langle Q,X \rangle\\
        \textrm{s.t.}~& 
         \mu^\top x \geq \rho,\\
        & e^\top x \leq 1,\\
        & e^\top y \geq n-\aleph, \\
        & \textrm{diag}(Z) = \textbf{0},\\
        &\textrm{diag}(Y) = y,\\
        &0 \leq x_i\leq u_i,  ~\forall i \in [n],\\
        & \bar{X} \succeq 0,
    \end{aligned}
\end{equation}
where $\bar{X} = \begin{pmatrix}1& x^\top &y^\top \\x &X &Z^\top\\ y& Z& Y \end{pmatrix}$. The positive semidefinite matrix in~\eqref{eq:SDP} is of dimension~$2n+1$.

\citet{madani-lowrank} proved that any polynomial optimization problem can be reformulated as a quadratically constrained quadratic problem (with certain sparsity properties) with a corresponding SDP relaxation having on optimal solution with rank at most two. Therefore, we can expect that the SDP relaxation~\eqref{eq:SDP} will yield strong bounds. In the next section we will evaluate the quality of the bounds of the SDP relaxation on problem~\eqref{eq:cardMIQP}.

\section{Numerical results}

We performed numerical experiments in order to evaluate the quality of the lower bounds that we obtain from the SDP relaxation~\eqref{eq:SDP}.
We used python and ran the tests on a ThinkPad-X1-Carbon-6th with 8~Intel(R) Core(TM) i7-8550U~CPU @~1.80GHz.

The model data $Q$, $\mu$, $\rho$, and upper bounds $u$ are from the instances in the paper of  \citet{frangioni2007sdp} and can be found at the website~\cite{WBfrangioni}.

We use Gurobi~9.1.2 \cite{gurobi} to solve the MIQP problem~\eqref{eq:cardMIQP}. We report the gap between the upper and lower bounds found by Gurobi after a time limit of 90~seconds.
The lower bounds from the SDP relaxation~\eqref{eq:SDP} are obtained by Mosek~\cite{mosek}.

The code can be downloaded from \href{https://github.com/shudianzhao/cardinality-SDP}{https://github.com/shudianzhao/cardinality-SDP}.

\input{average_results}


Table~\ref{tab:average} shows the overall performance on instances orl$n$-005- and  pard$n$-005- with $n\in \{200,300,400\}$ and $\aleph\in \{5,10,20\}$. We test 30~instances for each pair of $n$ and $\aleph$.

The optimality gap is measured as difference between the lower bound ($lb$) and the best found upper bounds ($ub$). We compute the relative gap as $(ub-lb)/ub$, where the upper bounds $ub$ are found by Gurobi and the lower bounds $lb$ are obtained from relaxations solved in Gurobi and the SDP relaxation~\eqref{eq:SDP} solved by Mosek.

When the time limit is reached, the gap for Gurobi solving~\eqref{eq:cardMIQP} is typically between 60\% and~80 \%. There is only one instance, namely orl300-05-i, that can be solved by Gurobi within 90~seconds for $\aleph\in\{5,20\}$.

In contrast to that, the lower bounds from the SDP relaxation~\eqref{eq:SDP} are exceptional. The gaps between the SDP lower bounds and the best upper bounds found by Gurobi are closed for most of the instances.

Actually, 96~\% of the SDP solutions have rank one, and hence the optimal solutions for problem~\eqref{eq:cardMIQP} can be built from it. The computation times for solving the SDP problems are less  than 40~seconds, even for large instances where $n=400$ and hence the size of the SDP matrix is $801$.

These results are another evidence of the strong modeling power of semidefinite programming, in particular when quadratic functions and binary variables are present.

\ifSpringer
\else
\bibliographystyle{plainnat}
\bibliography{mybib}

\appendix
\appendixpage

The detailed numerical results for all instances are given in Tables~\ref{tab:summary_200}, \ref{tab:summary_300} and~\ref{tab:summary_400}. In these tables we present the computational results of Gurobi and Mosek on solving~\eqref{eq:cardMIQP} and~\eqref{eq:SDP}, respectively. For the SDP solution, we also evaluate the rank of the positive semidefinite matrices by calculating their eigenvalues. 

\begin{landscape}

\begin{table}
	\centering
	\setlength{\tabcolsep}{2pt}
	\input{n200_compare}

\end{table}

\begin{table}
	\centering
	\setlength{\tabcolsep}{2pt}
	\input{n300_compare}

\end{table}

\begin{table}
	\centering
	\setlength{\tabcolsep}{2pt}
	\input{n400_compare}

\end{table}	
\end{landscape}

\fi

\end{document}

%% file: average_results.tex
\setlength{\tabcolsep}{10pt}
\begin{filecontents*}{cardinality_average_2.csv}
chi,n,min,Optimality Gap,max,min,gap,max,time
5,200,0.639196483680587,0.840537902789411,0.910456782417176,0,0.000243931536222,0.002327304215797,4.05312269528707
,300,0,0.860428837452379,0.93550277135987,0,0.000267043656075,0.004765123655825,13.2726219256719
,400,0.588070653728639,0.916234098564432,0.950891650810019,0,2.0541561358794E-05,0.000326185352248,30.107376964887
10,200,0.634042808545018,0.757798112614504,0.839839600769286,0,0.000151486995105,0.001548061378261,4.52696158885956
,300,0.137064737890971,0.802983190521068,0.882548157437901,0,0.000243438630544,0.004682181793085,12.2950302600861
,400,0.566214563241464,0.858283533760437,0.905745194180192,0,6.94022612621448E-05,0.000531426752893,27.9159217039744
20,200,0.416288090359734,0.593758671395664,0.705300763613384,0,0.000106017534831,0.000616612118192,4.24825749397278
,300,0,0.673748650816037,0.785960468683739,0,0.000113536894155,0.001284153230759,12.0523473421733
,400,0.332067962009432,0.748507534295108,0.827399804412574,0,5.64568915205761E-05,0.000345365272628,27.7429266532262
\end{filecontents*}
\pgfplotstabletypeset[  
col sep=comma, 
fixed zerofill,
fixed,
font=\footnotesize,
multicolumn names, 
display columns/0/.style={int detect,
	column name=$\aleph$, 
	assign column name/.style={
			/pgfplots/table/column name={\multicolumn{1}{c|}{##1}}}
},  
display columns/1/.style={int detect,
	column name=$n$, 
	assign column name/.style={
	/pgfplots/table/column name={\multicolumn{1}{c|}{##1}}}	
},
 display columns/3/.style={
column name= avg,
preproc/expr={100*##1},
}, 
display columns/2/.style={,
	column name=min,
	preproc/expr={100*##1},
},
display columns/4/.style={,
	column name=max,
	preproc/expr={100*##1},
		assign column name/.style={
		/pgfplots/table/column name={\multicolumn{1}{c|}{##1}}}
},
display columns/6/.style={,
	column name=avg ,
	preproc/expr={100*##1},
},
display columns/5/.style={,
	column name=min,
	preproc/expr={100*##1},
},
display columns/7/.style={,
	column name=max,
	preproc/expr={100*##1},
},
display columns/8/.style={,
	column name=time,
},
every head row/.style={
	fixed,precision=2,
	before row={
		\caption{Average results for Gurobi solving~\eqref{eq:cardMIQP} within a time limit of 90~seconds and Mosek for solving~\eqref{eq:SDP} on data with $n\in \{200,300,400\}$; 30 instances are tested for each pair of $n$ and $\aleph$ }\label{tab:average}\\
		\toprule  
		&  &  \multicolumn{3}{c|}{gap (MIQP~\eqref{eq:cardMIQP})} & \multicolumn{4}{c}{gap (SDP~\eqref{eq:SDP})} \\
	},
	after row={
		 & &  \% & \% &  \% & \% & \% & \%& (\si{\second})\\ 
		\endfirsthead
		\midrule} 
},
columns addvline/.style={
	every col no #1/.style={
	column type/.add={}{|}}
},
columns addvline/.list={0,1,4,},
every row no 3/.style ={before row={\midrule}},
every row no 6/.style ={before row={\midrule}},
every last row/.style={after row=\bottomrule},]{cardinality_average_2.csv} 

%% file: n200_compare.tex
\begin{filecontents*}{cardinality_summary_200.csv}
N 5,graph,bestObj,Optimality Gap,sdp,gap,time,rank,10-bestObj,Optimality Gap,sdp,gap,time,rank,20-bestObj,Optimality Gap,sdp,gap,time,rank
200,orl200-005-a,34.0388259175257,0.846640761884643,34.032761636783,0.000178189498915,3.7214469909668,2,18.5179940997104,0.736757014031888,18.5179941005468,-4.51636807891657E-11,3.75956416130066,,11.0645598104208,0.569367587124627,11.0614854399544,0.000277934684554,4.75709819793701,2
,orl200-005-b,42.0994131029935,0.853142670210509,42.0994131038523,-2.03996118975011E-11,3.59061288833618,,23.5451820184544,0.754456269542118,23.539908003809,0.000224045677859,4.61229801177979,2,13.7853602503118,0.588128574420945,13.7853602505228,-1.53105458799319E-11,4.01511287689209,
,orl200-005-c,71.4779052637416,0.835651215108002,71.4556042348272,0.000312096289063,3.40409588813782,2,38.1759412585974,0.724493636156069,38.1759412591755,-1.51435398061741E-11,3.31842613220215,,22.7630784878355,0.549201947223751,22.7629035775776,7.68400468996011E-06,4.53276991844177,2
,orl200-005-d,37.3657273735256,0.860852212080068,37.3657273752908,-4.72390953618937E-11,3.01419711112976,,20.3486380609861,0.759707809184343,20.3486380616046,-3.03916150861742E-11,3.73633408546448,,12.0071069427236,0.599129756801596,12.0071069428611,-1.14526432272889E-11,4.26180505752564,
,orl200-005-e,33.2523405000539,0.862141321307201,33.2523405012134,-3.48690597290659E-11,4.04484987258911,,17.9093978339912,0.759776251070686,17.9093978345528,-3.13583645200838E-11,3.74126791954041,,10.5426248795416,0.597794510627675,10.5426248799002,-3.40185365223263E-11,4.27822589874268,
,orl200-005-f,67.7271004112249,0.845936281194468,67.6993017944677,0.000410618957956,4.62415099143982,2,37.219221235591,0.749137001796003,37.2115097945836,0.000207232682845,4.3235490322113,2,21.6111968661019,0.576871315615575,21.6105842272798,2.83490171172223E-05,4.23996710777283,2
,orl200-005-g,71.8480148379141,0.844016787374393,71.8480148426885,-6.64514761927004E-11,3.83506393432617,,40.7890620350396,0.754174976541042,40.7782102660297,0.000266116853561,4.56311917304993,2,23.9903130229817,0.590634036820069,23.9902266030971,3.60229547002993E-06,4.68522095680237,2
,orl200-005-h,89.1108332190979,0.855764219125584,89.1108332211804,-2.33693095967004E-11,4.05064487457275,,48.846807796905,0.763840373974334,48.8468077981725,-2.59487963937695E-11,3.62124514579773,,27.9281034424949,0.595940355502785,27.9274845023589,2.21624019149423E-05,4.25100493431091,2
,orl200-005-i,10.4367550576164,0.867131909559609,10.4367550577715,-1.48610195618131E-11,4.39112091064453,,5.81622737284971,0.770696675236795,5.81523455021072,0.000170727875276,4.54070997238159,2,3.35347128997286,0.608208226122463,3.35255492079507,0.000273334576001,4.6289529800415,2
,orl200-005-j,38.2441287916768,0.867736541064822,38.2441287934385,-4.6062871596059E-11,4.14484190940857,,21.0036472162785,0.77055664675114,21.0036472170722,-3.77873744939517E-11,3.66861891746521,,12.0013839614656,0.605241522062419,12.0013839621734,-5.89811233515352E-11,3.94467091560364,
,orl200-05-a,17.9474877002351,0.820087607876391,17.9474877007834,-3.05471231356536E-11,3.78273415565491,,10.1215705316561,0.703038997483568,10.1215705320653,-4.04281362420372E-11,3.84577703475952,,6.27663682638043,0.523023922718403,6.27571977429222,0.000146126997571,4.57679009437561,2
,orl200-05-b,10.3379390621707,0.85244088918193,10.3355083028737,0.000235185268664,4.88529205322266,2,5.66081322597388,0.741068093481323,5.6608132261983,-3.964378126221E-11,4.31163787841797,,3.35165058299502,0.569556431804374,3.35136085438117,8.64510347979448E-05,4.58993911743164,2
,orl200-05-c,12.0100819949175,0.854102189406481,12.0092514991077,6.91546687880606E-05,4.62750697135925,,6.60742161790232,0.744441407741128,6.60728877952019,2.01048246194987E-05,4.75598096847534,2,3.88424306846702,0.571280923947347,3.88321360197316,0.000265106841749,4.40048813819885,2
,orl200-05-d,97.570970474482,0.756152142355248,97.3444203945136,0.002327304215797,4.37479496002197,2,53.5305323736213,0.709706227030929,53.4477919112102,0.001548061378261,3.95942306518555,2,30.8354668231568,0.516563748109211,30.8164650173873,0.000616612118192,4.35039496421814,2
,orl200-05-e,108.59128384973,0.639196483680587,108.426749349648,0.001517471482533,4.12604999542236,2,59.6835643651486,0.634042808545018,59.6419097125866,0.000698412454644,4.35385990142822,2,34.7837158388239,0.416288090359734,34.7810604246014,7.63465572962394E-05,4.50139999389648,2
,orl200-05-f,81.7713656970768,0.790428760965541,81.7129395417655,0.000715017176458,4.19312405586243,2,45.2396383610774,0.712758427227616,45.2319895749151,0.000169101254094,4.49461483955383,2,26.9845605293284,0.533202135703283,26.971242050912,0.000493802932445,4.20925116539001,2
,orl200-05-g,5.75919118940951,0.855937231400445,5.7591889120959,3.95422627888911E-07,4.69707703590393,2,3.23221176894211,0.750669633839865,3.2303344068279,0.000581166491691,4.46881318092346,2,1.88174633552602,0.577718690995508,1.88124443138167,0.000266793690378,4.52176380157471,2
,orl200-05-h,21.2001625367448,0.820888169909113,21.192373325968,0.000367547827564,4.45833802223206,2,11.9016058827153,0.702799971801561,11.9016058829109,-1.64362478601442E-11,4.2353310585022,,7.21242294693248,0.517689499144935,7.21242294739099,-6.35715954436746E-11,3.99059987068176,
,orl200-05-i,20.6265071257704,0.841314760686896,20.6160499973424,0.000507232395602,4.48742198944092,,11.7535798657525,0.734714976409584,11.7491790062404,0.000374567406782,4.42264580726624,2,7.01626496702451,0.561023395603792,7.01587437775835,5.56722149130514E-05,4.69561696052551,2
,orl200-05-j,64.1920680767598,0.845591557716734,64.1502027837225,0.000652613572843,4.20066499710083,2,35.7999109760228,0.74541782917912,35.7929184026644,0.000195361922708,4.35875916481018,2,21.1474766274497,0.575345594719177,21.1362265519756,0.000532265087453,4.17830300331116,2
,pard200-a,141.034902960883,0.910456782417176,141.034902963377,-1.76858181124809E-11,5.60545516014099,,74.6274341041905,0.839839600769286,74.6274089164447,3.37513336300862E-07,4.5387179851532,2,40.1248822058333,0.705300763613384,40.1248704352241,2.93349462197288E-07,4.44686508178711,2
,pard200-b,381.186332750074,0.782277230704329,381.186332767477,-4.56563562916962E-11,3.52814292907715,,207.023298355397,0.751075418475016,207.023298356549,-5.56508577834863E-12,6.32732105255127,,115.155703905809,0.578836601715201,115.155703909645,-3.33066683976296E-11,3.72623491287231,2
,pard200-c,356.04434214972,0.763117169232084,356.044342161179,-3.21823000308147E-11,3.26545310020447,,194.895233679078,0.734421216245924,194.895233871609,-9.87870614658838E-10,6.27848315238953,,109.985811494671,0.559076244797155,109.982947402899,2.60412349341389E-05,4.22175598144531,
,pard200-d,342.407330990136,0.800705352807327,342.407331008329,-5.31313556359064E-11,3.32660818099976,,184.068059630161,0.777904372863715,184.068059714243,-4.56799912590327E-10,8.14832901954651,,103.317393679855,0.619633618009105,103.31739368545,-5.41623006308582E-11,3.75582599639893,2
,pard200-e,101.799614356857,0.890499771854507,101.799614360253,-3.33568273047681E-11,3.94435501098633,,55.8364232334696,0.814091626283889,55.8314563424339,8.896223314046E-05,4.4853630065918,2,32.0353895457431,0.681609056541003,32.0353271580806,1.94746450341556E-06,4.74158000946045,
,pard200-f,25.0854950458396,0.900015570665047,25.0848649160383,2.51199200536743E-05,4.61118197441101,2,13.6856078219814,0.825230578320628,13.6856078227679,-5.74667775063193E-11,3.79059100151062,,7.70844796932461,0.695411634416241,7.70844796977379,-5.82710504286566E-11,3.92002105712891,
,pard200-g,324.538489861692,0.857381968339513,324.538489871954,-3.16174816728434E-11,3.4641330242157,,177.763432940405,0.805404608513259,177.763432979028,-2.17269188005148E-10,5.85192489624023,,100.422462340083,0.667808261826991,100.422462344825,-4.72252397176061E-11,3.5932879447937,
,pard200-h,55.5869513477572,0.894883604637606,55.5869513505268,-4.98234224377445E-11,3.57020401954651,,30.2215251187201,0.813200107868605,30.2215126269247,4.13341169409947E-07,4.89555311203003,2,17.0629123565981,0.673885211915204,17.06291235696,-2.12089538514427E-11,3.859375,
,pard200-i,130.410469628629,0.902796926148672,130.410469631184,-1.95863034136305E-11,3.63505601882935,,70.0280458242337,0.828993454425721,70.0280458256005,-1.95180973743976E-11,3.84371304512024,,39.4232818910446,0.700362498133849,39.423281892243,-3.03986764030293E-11,3.70194101333618,
,pard200-j,71.3026020142383,0.898848994787402,71.30260201733,-4.33612562204498E-11,3.98906183242798,,38.4227177764571,0.821527367644947,38.422717777612,-3.00594535602825E-11,4.55687594413757,,21.7156799643294,0.6886259854741,21.7156799649941,-3.0612259025753E-11,3.87146186828613,
\end{filecontents*}
\pgfplotstabletypeset[  
col sep=comma, 
fixed zerofill,
fixed,
font=\footnotesize,
multicolumn names, 
display columns/0/.style={int detect,
	column name=$n$, 
	},  
display columns/1/.style={string type,
	column name=name, 
 	assign column name/.style={
	/pgfplots/table/column name={\multicolumn{1}{c|}{##1}}}
},  
display columns/2/.style={
	column name=ub
},
display columns/3/.style={
	column name=gap,
	preproc/expr={100*##1},
		assign column name/.style={
		/pgfplots/table/column name={\multicolumn{1}{c|}{##1}}}
},
display columns/4/.style={,
	column name=$lb_{SDP}$,
},
display columns/6/.style={,
	column name=time,
},
display columns/5/.style={,
	column name=gap,
	preproc/expr={100*##1},
},
display columns/7/.style={int detect,
	empty cells with={1},
	column name=rank,
	assign column name/.style={
		/pgfplots/table/column name={\multicolumn{1}{c|}{##1}}}
},
display columns/8/.style={
	column name=ub,
},
display columns/9/.style={
	column name=gap,
	preproc/expr={100*##1},
	assign column name/.style={
		/pgfplots/table/column name={\multicolumn{1}{c|}{##1}}}
},
display columns/10/.style={,
	column name=$lb_{SDP}$,
},
display columns/12/.style={,
	column name=time,
},
display columns/11/.style={
	preproc/expr={100*##1},
	column name=gap,
},
display columns/13/.style={int detect,
	empty cells with={1},
	column name=rank,
	assign column name/.style={
		/pgfplots/table/column name={\multicolumn{1}{c|}{##1}}}
},
display columns/14/.style={
	column name=ub,
},
display columns/15/.style={
	preproc/expr={100*##1},
	column name=gap,	assign column name/.style={
		/pgfplots/table/column name={\multicolumn{1}{c|}{##1}}}
},
display columns/16/.style={,
	column name=$lb_{SDP}$,
},
display columns/18/.style={,
	column name=time,
},
display columns/17/.style={
	preproc/expr={100*##1},
	column name=gap,
},
display columns/19/.style={int detect,
	empty cells with={1},
	column name=rank,
},
every head row/.style={
	fixed,precision=2,
	before row={
		\caption{Detailed comparisons of experiments on instances with $n=200$ }\label{tab:summary_200}\\
		\toprule  
		&  &  \multicolumn{6}{c|}{$\aleph = 5$}& \multicolumn{6}{c|}{$\aleph = 10$}& \multicolumn{6}{c}{$\aleph = 20$} 
		\\
		& & \multicolumn{2}{c|}{Gurobi} &  \multicolumn{4}{c|}{Mosek} &\multicolumn{2}{c|}{Gurobi} & \multicolumn{4}{c|}{Mosek}& \multicolumn{2}{c|}{Gurobi} &  \multicolumn{4}{c}{Mosek}\\
	},
after row={
& & & \% &  &  \% & & & &  \% & &  \% & & & &\% & & \% & &\\ 
		\endfirsthead
		\midrule} 
},
columns addvline/.style={
	every col no #1/.style={
	column type/.add={}{|}}
},
columns addvline/.list={1,3,7,9,13,15},
every row no 30/.style ={before row={\midrule}},
every row no 60/.style ={before row={\midrule}},
every last row/.style={after row=\bottomrule},
]{cardinality_summary_200.csv} 

%% file: n300_compare.tex
\begin{filecontents*}{cardinality_summary_300.csv}
N 5,graph,bestObj,Optimality Gap,sdp,gap,time,rank,10-bestObj,Optimality Gap,sdp,gap,time,rank,20-bestObj,Optimality Gap,sdp,gap,time,rank
300,orl300-005-a,207.755319416486,0.857034118695663,207.755319427962,-5.52393073445949E-11,8.30320596694946,,110.951838603905,0.786681967508689,110.951838606673,-2.49470334328196E-11,10.5142431259155,,62.2487198720169,0.633192535207385,62.2470070047284,2.75172634143438E-05,11.9344830513,2
,orl300-005-b,99.3675146918572,0.906195998831587,99.3675146975513,-5.73033405308239E-11,7.87116599082947,,54.1593473831177,0.830524288129344,54.1448155610901,0.00026838806037,15.1257450580597,2,30.1456304146111,0.698732459853248,30.1426002090953,0.00010052900197,12.2599258422852,2
,orl300-005-c,43.6012264880366,0.906968193387933,43.6012264910372,-6.88191800625069E-11,10.8694288730621,,23.6215547591541,0.831264500204319,23.6208071974734,3.16484392134165E-05,13.9572811126709,2,13.2013369238602,0.700850894692361,13.1994776569625,0.000140859126853,13.3490500450134,2
,orl300-005-d,168.984395115847,0.894052187360751,168.984395124192,-4.93834960812507E-11,15.6845350265503,,95.1631537199637,0.828935564049292,95.1383767575055,0.000260430788318,12.9444050788879,2,53.9693844554608,0.703493879510001,53.9580535199703,0.000209995260232,12.2186589241028,2
,orl300-005-e,95.1081628386338,0.910544002474004,95.1081628426271,-4.19862590630054E-11,11.7387568950653,,50.4083890628337,0.835362631176635,50.4083890654374,-5.16522297927824E-11,11.4187440872192,,28.0558095793558,0.705516468424421,28.0558095809455,-5.66615837617708E-11,10.0629079341888,
,orl300-005-f,62.5444397404906,0.9077827509898,62.5444397419164,-2.27967599773293E-11,12.9132099151611,,34.4251122314057,0.836420818673939,34.422073319355,8.82838178452984E-05,13.5738801956177,2,19.6745576706106,0.716926200529825,19.6737174317832,4.27086965274899E-05,13.3905701637268,2
,orl300-005-g,39.0960224318252,0.90551398731498,39.0959150911527,2.74557258236371E-06,13.5553550720215,2,20.8871407255377,0.827639568338684,20.8871407268724,-6.3900449714945E-11,11.4362308979034,,11.8429314570526,0.697943994307141,11.8423559405309,4.85981442054195E-05,12.4584159851074,2
,orl300-005-h,45.1560848914567,0.902692022450505,45.1509011039019,0.000114810279041,12.7488329410553,2,24.6543595844119,0.825812632066896,24.6543595847665,-1.43822715378555E-11,11.9248838424683,,13.9175620059253,0.694677305561332,13.9175620069142,-7.10589347415545E-11,10.2985649108887,
,orl300-005-i,38.4357362909012,0.91188603958593,38.4146908819935,0.000547847930687,11.9961569309235,2,20.4421013673176,0.837387817785423,20.4420703477445,1.51743794071049E-06,13.7762830257416,2,11.3431238490572,0.709228285603249,11.3405427376341,0.000227600343548,13.0265820026398,2
,orl300-005-j,15.2652047941694,0.908586346561918,15.2627367343739,0.000161704931327,12.3241589069366,2,8.11336283136317,0.831128682388811,8.1133628316274,-3.25667642279509E-11,12.9816739559174,,4.57716491978872,0.702724376159941,4.57708479333581,1.75060014257673E-05,12.9594731330872,2
,orl300-05-a,157.743559047635,0.704053211406234,157.743559054128,-4.11599192045138E-11,10.0682990550995,,86.8062506627402,0.666193309351159,86.7718874918985,0.000396017325829,12.9490959644318,2,50.005670486034,0.460289747786794,49.9931781156366,0.000249881501203,13.1186017990112,2
,orl300-05-b,111.434778027415,0.884170420246505,111.32696947395,0.000968395654475,11.5008389949799,2,59.8032215553703,0.797781439853498,59.7951753525796,0.000134562742618,13.9678790569305,2,34.6035081872694,0.656027408565458,34.6030817769964,1.23228987447258E-05,12.5786719322205,2
,orl300-05-c,53.1759333716062,0.893546098023628,53.1759333729712,-2.56691501309573E-11,9.72164511680603,,29.0113552523826,0.807603168844973,29.0113552532102,-2.85277786318082E-11,12.4731781482697,,16.6474322409645,0.666893950192128,16.6444345395229,0.000180102329969,12.2823870182037,2
,orl300-05-d,113.849381853311,0.898525328128921,113.773230613379,0.000669324756995,11.383278131485,2,61.1639120247735,0.81891052156514,61.1452903202467,0.000304548468563,12.8067500591278,2,33.9828650527345,0.677613367403146,33.9808722952916,5.86435046637576E-05,12.209566116333,
,orl300-05-e,33.2624294071994,0.895220076179492,33.2566614878294,0.000173436512023,12.2615029811859,2,18.1628030048517,0.812483677314255,18.1528648892206,0.000547468165036,13.7409958839417,2,10.2393910038508,0.669982577140106,10.2393737229058,1.68769550078508E-06,13.1378200054169,2
,orl300-05-f,32.6320610636125,0.900790182697735,32.6320610673975,-1.15989918359832E-10,9.81165313720703,,17.6232411527796,0.818988796005442,17.6212047807511,0.000115563723016,12.853492975235,2,9.91582924333965,0.679897519565759,9.91277037842957,0.00030857820703,13.0434069633484,2
,orl300-05-g,173.625165837861,0.553906842655352,173.588477678526,0.000211351351343,12.0542159080505,2,93.6078627851434,0.603789825295122,93.6078627877466,-2.78099462409664E-11,12.7964270114899,,53.9971450311749,0.369234183356158,53.9830005473277,0.000262017370354,13.103856086731,2
,orl300-05-h,33.8606133495919,0.898844737985438,33.8606133503308,-2.18199398768673E-11,11.1105570793152,,18.5259464831596,0.819495751990131,18.5182210752525,0.000417178727682,13.3870739936829,2,10.4216810940674,0.681847346614318,10.4195155980589,0.000207830775636,12.9045751094818,2
,orl300-05-i,179.358392046119,0,178.507780398991,0.004765123655825,14.2694528102875,2,99.0417370127381,0.137064737890971,98.5801667508181,0.004682181793085,12.3076000213623,2,58.2722181920433,0,58.197483705321,0.001284153230759,12.5353419780731,2
,orl300-05-j,139.087885136783,0.861992028844517,139.032748773318,0.000396571052157,15.1153290271759,2,75.1761410671133,0.784019040594601,75.1761410684183,-1.73588130368999E-11,9.87188601493835,,42.9924930233802,0.632184474660528,42.9924930244727,-2.54114616747955E-11,10.346174955368,
,pard300-a,64.7869844772603,0.931644917042448,64.7869844785272,-1.95545987701501E-11,11.8460412025452,,34.4054571385364,0.87362397022677,34.4054571391869,-1.89067323398566E-11,10.7901251316071,,19.1266180178758,0.774052417380469,19.1266180183905,-2.69101099693282E-11,10.8481130599976,
,pard300-b,232.074900885452,0.926504239143533,232.074900971922,-3.72592328321323E-10,16.9974918365479,,125.257281418567,0.869452056428298,125.257281420866,-1.83522094097449E-11,12.1936869621277,,69.8051515932125,0.766491413935942,69.8051515945603,-1.93078283208573E-11,10.124410867691,
,pard300-c,375.475736133038,0.93550277135987,375.475736128036,1.33231114437168E-11,21.1155610084534,2,196.500773233371,0.881636723099447,196.500773244872,-5.85254934045607E-11,9.28862905502319,,106.365673152492,0.783078131452184,106.36371006931,1.845632481898E-05,12.1363170146942,2
,pard300-d,195.903151589439,0.934653211665108,195.903151601099,-5.95212690681891E-11,11.4335339069366,,106.346602289922,0.882548157437901,106.346602292757,-2.66537872668021E-11,10.428297996521,,57.7595373208373,0.785960468683739,57.7595373217293,-1.54423725943543E-11,11.3786880970001,
,pard300-e,113.962621205035,0.932027291343376,113.962621227032,-1.93016651853692E-10,19.7437679767609,,60.5888823269118,0.874582828601819,60.5888823310241,-6.78719798300024E-11,10.2595269680023,,33.4992867468304,0.775482982140193,33.499286747331,-1.49446226881753E-11,11.0610921382904,
,pard300-f,53.4660170775881,0.934194532727602,53.4660170855529,-1.48970593992344E-10,20.4242401123047,,28.4330965638906,0.877722735883113,28.4330965643899,-1.75604644146974E-11,11.2057371139526,,15.7397540315126,0.781100465900028,15.7397540321174,-3.84248691709753E-11,10.9727880954742,
,pard300-g,467.22981466513,0.920420109942592,467.229814748131,-1.77645030984868E-10,15.7534260749817,,253.059097019805,0.865054399583329,253.059097135321,-4.56477350980124E-10,16.472892999649,,139.345571715811,0.759481869469481,139.345340187701,1.66154182176833E-06,12.4784359931946,2
,pard300-h,81.6854803459144,0.93244988900016,81.6854803503701,-5.45463275307142E-11,10.8985350131989,,44.0044115084691,0.876986576369702,44.0019751005758,5.53704211616539E-05,12.5847859382629,2,24.2292725677348,0.777135087602936,24.2292725692836,-6.39221734897706E-11,10.4983441829681,
,pard300-i,310.311283657211,0.928935048759135,310.311283789613,-4.26676686094093E-10,16.9578638076782,,165.493903813506,0.872135191633533,165.49390381974,-3.76683097928959E-11,9.97083902359009,,92.4125963068454,0.772551813037985,92.4125963085993,-1.89787468083244E-11,11.6967349052429,
,pard300-j,150.822371052207,0.934228538766653,150.822371063492,-7.48250700285172E-11,17.7066180706024,,80.0095367772391,0.878264337340802,80.0095367797342,-3.11846676914342E-11,10.8486371040344,,43.8016382528166,0.779867899744839,43.8013991866618,5.45795703499049E-06,13.1564619541168,2
\end{filecontents*}
\pgfplotstabletypeset[  
col sep=comma, 
fixed zerofill,
fixed,
font=\footnotesize,
multicolumn names, 
display columns/0/.style={int detect,
	column name=$n$, 
	},  
display columns/1/.style={string type,
	column name=name, 
 	assign column name/.style={
	/pgfplots/table/column name={\multicolumn{1}{c|}{##1}}}
},  
display columns/2/.style={
	column name=ub
},
display columns/3/.style={
	column name=gap,
	preproc/expr={100*##1},	assign column name/.style={
		/pgfplots/table/column name={\multicolumn{1}{c|}{##1}}}
},
display columns/4/.style={,
	column name=$lb_{SDP}$,
},
display columns/6/.style={,
	column name=time,
},
display columns/5/.style={,
	column name=gap,
	preproc/expr={100*##1},
},
display columns/7/.style={int detect,
	empty cells with={1},
	column name=rank,
	assign column name/.style={
		/pgfplots/table/column name={\multicolumn{1}{c|}{##1}}}
},
display columns/8/.style={
	column name=ub,
},
display columns/9/.style={
	column name=gap,
	preproc/expr={100*##1},	assign column name/.style={
		/pgfplots/table/column name={\multicolumn{1}{c|}{##1}}}
},
display columns/10/.style={,
	column name=$lb_{SDP}$,
},
display columns/12/.style={,
	column name=time,
},
display columns/11/.style={
	preproc/expr={100*##1},
	column name=gap,
},
display columns/13/.style={int detect,
	empty cells with={1},
	column name=rank,
	assign column name/.style={
		/pgfplots/table/column name={\multicolumn{1}{c|}{##1}}}
},
display columns/14/.style={
	column name=ub,
},
display columns/15/.style={
	preproc/expr={100*##1},
	column name=gap,	assign column name/.style={
		/pgfplots/table/column name={\multicolumn{1}{c|}{##1}}}
},
display columns/16/.style={,
	column name=$lb_{SDP}$,
},
display columns/18/.style={,
	column name=time,
},
display columns/17/.style={
	preproc/expr={100*##1},
	column name=gap,
},
display columns/19/.style={int detect,
	empty cells with={1},
	column name=rank,
},
every head row/.style={
	fixed,precision=2,
	before row={
		\caption{Detailed comparisons of experiments on instances with $n=300$}\label{tab:summary_300}\\
		\toprule  
		&  &  \multicolumn{6}{c|}{$\aleph = 5$}& \multicolumn{6}{c|}{$\aleph = 10$}& \multicolumn{6}{c}{$\aleph = 20$} \\
	& & \multicolumn{2}{c|}{Gurobi} &  \multicolumn{4}{c|}{Mosek} &\multicolumn{2}{c|}{Gurobi} & \multicolumn{4}{c|}{Mosek}& \multicolumn{2}{c|}{Gurobi} &  \multicolumn{4}{c}{Mosek}\\
	},
after row={
& & & \% &  &  \% & & & &  \% & &  \% & & & &\% & & \% & &\\ 
		\endfirsthead
		\midrule} 
},
columns addvline/.style={
	every col no #1/.style={
	column type/.add={}{|}}
},
columns addvline/.list={1,3,7,9,13,15},
every row no 30/.style ={before row={\midrule}},
every row no 60/.style ={before row={\midrule}},
every last row/.style={after row=\bottomrule},
]{cardinality_summary_300.csv} 

%% file: n400_compare.tex
\begin{filecontents*}{cardinality_summary_400.csv}
N 5,graph,bestObj,Optimality Gap,sdp,gap,time,rank,10-bestObj,Optimality Gap,sdp,gap,time,rank,20-bestObj,Optimality Gap,sdp,gap,time,rank
400,orl400-005-a,202.68182803638,0.933991308218229,202.681828062345,-1.28108143245273E-10,31.9379918575287,,106.648499587429,0.875419153104297,106.641614245703,6.45652428880601E-05,28.5606200695038,,58.1089947054396,0.771489927239237,58.1089947083234,-4.9627457098574E-11,22.8488168716431,
,orl400-005-b,247.846926142699,0.918501936441544,247.846926160721,-7.27151821114042E-11,37.9740560054779,,136.416860445753,0.85910906347194,136.408224359612,6.33105971526519E-05,27.6458859443665,,75.7092803367934,0.747589358124321,75.709280339964,-4.18791543743405E-11,22.3083350658417,
,orl400-005-c,114.928696104468,0.936173370000596,114.92869610571,-1.08081867518385E-11,27.8868899345398,,60.9357228691328,0.880273851990215,60.9178037068846,0.000294153123682,25.9828259944916,,32.9080215360174,0.778446766557263,32.9079488861624,2.20766888921667E-06,29.3602070808411,2
,orl400-005-d,96.0636389047301,0.931372439684623,96.0636389067211,-2.07267065541866E-11,26.147696018219,,51.1828009642522,0.871098975367134,51.1816263043693,2.29508119934025E-05,29.809632062912,,28.3417815662562,0.767388893403297,28.3417815666895,-1.52886198355355E-11,25.8489091396332,
,orl400-005-e,133.89078543206,0.9314952023183,133.89078543543,-2.517542849618E-11,25.359475851059,,71.2101217720171,0.871494660321599,71.2033316812905,9.53619804893523E-05,26.8182890415192,,39.5009944462788,0.768517725559685,39.5009944487715,-6.31055791330454E-11,25.3703300952911,
,orl400-005-f,190.12886767311,0.936673389037486,190.089092393559,0.000209245459853,36.5676889419556,2,100.148493837128,0.879703441792675,100.135604520395,0.000128718619058,27.6162719726562,2,54.3295078842943,0.778683131555032,54.3295078851661,-1.60460129654414E-11,27.2882299423218,
,orl400-005-g,115.351919739384,0.927223826733437,115.351919741889,-2.17137260704266E-11,27.3976650238037,,61.9029229261919,0.864578168142745,61.9029229286589,-3.98519895968507E-11,22.9083340167999,,34.3586962649697,0.755979453082919,34.3586330496438,1.83986731400472E-06,29.949550151825,2
,orl400-005-h,41.4796618000054,0.93431372535123,41.4778891310477,4.27376849419906E-05,36.2644739151001,2,21.8224338544043,0.87611620964328,21.8224338547221,-1.45645015239249E-11,24.0654900074005,2,11.854869549747,0.772627035487679,11.8540171501184,7.19080812712772E-05,29.0552401542664,2
,orl400-005-i,188.037265364466,0.931329371948773,188.037265368766,-2.2869802741723E-11,20.9794352054596,,98.8960761896159,0.870312651313569,98.896076191668,-2.07508264998334E-11,22.9023711681366,,54.3739244507268,0.764371856841791,54.3716739406779,4.13912224098479E-05,28.2515108585358,2
,orl400-005-j,245.54733395398,0.903800482064795,245.547333988809,-1.41841553811728E-10,40.8693740367889,,131.884421141865,0.835498911940733,131.864536331399,0.000150797257694,32.9654719829559,,72.9120456772798,0.706168685084186,72.9114867922283,7.66525380419884E-06,28.1922390460968,2
,orl400-05-a,134.528333079142,0.918516152500598,134.528333083467,-3.21499087473883E-11,19.9424860477448,,72.6165546957802,0.850439626898382,72.5959055824214,0.000284439090513,27.7293920516968,,40.5483799229703,0.732637989501416,40.5468993609081,3.65148034880352E-05,26.9439759254456,2
,orl400-05-b,105.792937047607,0.92446975500791,105.788908138213,3.80844217467872E-05,28.1261079311371,,56.8560072624919,0.860230096495399,56.8560072630539,-9.88568574060862E-12,21.5071890354156,2,31.6604160664486,0.749179643001251,31.655598923732,0.000152173482112,25.9999778270721,2
,orl400-05-c,122.822250002724,0.919211406250053,122.822250011563,-7.19715349243098E-11,31.1371719837189,,65.8515845999589,0.85039268270337,65.8227617206713,0.000437886204318,26.7325961589813,,36.0359184826702,0.727459286028001,36.0330748642812,7.89168950940216E-05,31.32608294487,2
,orl400-05-d,176.900032184292,0.913799769012068,176.900032194115,-5.55254704715995E-11,21.4861629009247,,96.2254748594679,0.848272546529958,96.2246607421815,8.46058879399829E-06,29.9651510715485,,53.4852140211714,0.729491900729334,53.468501077613,0.000312575501868,25.840854883194,2
,orl400-05-e,171.06682740235,0.913997522631446,171.066827411311,-5.23810172683319E-11,20.5537519454956,,92.7968753749975,0.847084085083881,92.7968753812958,-6.78719798195477E-11,24.5525560379028,,52.2781169188464,0.730411428327549,52.2721069128785,0.000114975391711,28.259654045105,2
,orl400-05-f,76.8491004634466,0.912683940000742,76.8491004658131,-3.07939974419245E-11,21.3261818885803,,42.9458277748584,0.844926284468281,42.9458277754145,-1.29478592672387E-11,24.7119798660278,,24.8280343151218,0.732061135664371,24.819462534677,0.000345365272628,26.5574429035187,2
,orl400-05-g,174.020883139703,0.915002213184239,174.020883142265,-1.47250884736686E-11,20.2432770729065,,92.5678879207893,0.845789392992185,92.5678879224607,-1.80554279828619E-11,22.8329000473022,,51.4362985955427,0.72613317528069,51.428958548569,0.000142722061284,27.914039850235,2
,orl400-05-h,82.3589541278141,0.922652396407324,82.3589541314946,-4.46890123335877E-11,27.7380270957947,,44.0963960179704,0.855635894336416,44.0963960188923,-2.09058435486408E-11,22.9640769958496,,24.6319745098388,0.741790989126593,24.6319745111846,-5.46364428907366E-11,24.2277300357819,
,orl400-05-i,34.0140639533218,0.918940870161588,34.0029726817001,0.000326185352248,27.7498829364777,2,18.3551267326579,0.850057733345134,18.355126733233,-3.13358602648806E-11,23.4063310623169,2,10.3751931142541,0.73498229423911,10.3745377464272,6.31707978665801E-05,29.222470998764,2
,orl400-05-j,38.6634387512451,0.922198544540897,38.66343875355,-5.96148569290372E-11,23.9975850582123,,20.6166732000905,0.854175859225327,20.6166732003962,-1.48290279398532E-11,24.8312571048737,,11.5392704756237,0.739803106620119,11.5392397558781,2.66219839549628E-06,30.3876051902771,2
,pard400-a,318.7582064777,0.950320741382253,318.758206837706,-1.12940212855688E-09,37.1773250102997,,165.171834259348,0.904201534245642,165.171834273706,-8.69267563656385E-11,33.9217939376831,,88.6342176406543,0.82159945079675,88.6342176434543,-3.15902705505265E-11,23.6606359481812,
,pard400-b,861.624921717672,0.588070653728639,861.624921750398,-3.79818638645415E-11,28.2475800514221,,451.992079412239,0.566214563241464,451.752006310413,0.000531426752893,48.8205409049988,,254.620846253218,0.332067962009432,254.540388229956,0.000316091382674,33.9619691371918,2
,pard400-c,716.54230935799,0.919664371945901,716.542309405102,-6.57491571776511E-11,21.5555601119995,,376.214019018965,0.859732097436528,376.214019035538,-4.40520080257758E-11,26.129958152771,,203.077513072143,0.745776479338993,203.077513080787,-4.25663277048924E-11,28.1187810897827,
,pard400-d,373.442185038016,0.950262252811518,373.442186042515,-2.68983937188132E-09,49.0467510223389,,195.469654831022,0.905745194180192,195.469654887265,-2.87731093824893E-10,35.9624319076538,,106.554279356131,0.827273794947458,106.554279362586,-6.05823549716254E-11,24.6038341522217,
,pard400-e,290.890241143617,0.950891650810019,290.890241336904,-6.64469097617149E-10,36.8249099254608,,150.894689978437,0.905350805640906,150.894689988161,-6.44459041321769E-11,22.3111419677734,,82.4966025820671,0.827075063661344,82.4966025834908,-1.7258028160085E-11,31.3892059326172,
,pard400-f,802.745236637509,0.86849179299083,802.745236697302,-7.44861036759478E-11,43.9114720821381,,431.084949427457,0.804666545899581,431.084950047359,-1.43800423123688E-09,37.5734081268311,,233.703311503138,0.655456724032184,233.703311531308,-1.20536744031794E-10,40.956817150116,
,pard400-g,203.549395674823,0.948296307279089,203.54939572052,-2.24503127919851E-10,36.8303010463715,,108.918815006598,0.903445621726284,108.918815013502,-6.33845982288409E-11,23.6591517925262,,60.3369601066553,0.825911297546654,60.336960110719,-6.7350691962309E-11,24.9522290229797,
,pard400-h,438.276112216618,0.947931304576141,438.276112318397,-2.32224701537622E-10,34.9966239929199,,233.890285845939,0.903905520470965,233.890285867397,-9.1740679489968E-11,38.16490483284,,130.030343244774,0.827399804412574,130.029884559994,3.52753354843993E-06,29.9958829879761,2
,pard400-i,61.9605564652898,0.946450348141329,61.9605564677384,-3.95199723985499E-11,23.6022810935974,,33.1109789934966,0.899838340120655,33.1109789950111,-4.57416334083714E-11,23.8708119392395,,18.0507126827073,0.816412584169562,18.0507126838102,-6.1098554907081E-11,24.0805070400238,
,pard400-j,98.6842327659967,0.950295911771359,98.6842327772508,-1.14040817662553E-10,37.343122959137,,51.4822321966174,0.904796500684367,51.4822321996334,-5.85852140638384E-11,28.5548858642578,,27.378589810639,0.821039086484453,27.378589811383,-2.71745987763663E-11,25.4147341251373,
\end{filecontents*}
\pgfplotstabletypeset[  
col sep=comma, 
fixed zerofill,
fixed,
font=\footnotesize,
multicolumn names, 
display columns/0/.style={int detect,
	column name=$n$, 
	},  
display columns/1/.style={string type,
	column name=name, 
 	assign column name/.style={
	/pgfplots/table/column name={\multicolumn{1}{c|}{##1}}}
},  
display columns/2/.style={
	column name=ub
},
display columns/3/.style={
	column name=gap,
	preproc/expr={100*##1},
		assign column name/.style={
		/pgfplots/table/column name={\multicolumn{1}{c|}{##1}}}
},
display columns/4/.style={,
	column name=$lb_{SDP}$,
},
display columns/6/.style={,
	column name=time,
},
display columns/5/.style={,
	column name=gap,
	preproc/expr={100*##1},
},
display columns/7/.style={int detect,
	empty cells with={1},
	column name=rank,
	assign column name/.style={
		/pgfplots/table/column name={\multicolumn{1}{c|}{##1}}}
},
display columns/8/.style={
	column name=ub,
},
display columns/9/.style={
	column name=gap,
	preproc/expr={100*##1},
		assign column name/.style={
		/pgfplots/table/column name={\multicolumn{1}{c|}{##1}}}
},
display columns/10/.style={,
	column name=$lb_{SDP}$,
},
display columns/12/.style={,
	column name=time,
},
display columns/11/.style={
	preproc/expr={100*##1},
	column name=gap,
},
display columns/13/.style={int detect,
	empty cells with={1},
	column name=rank,
	assign column name/.style={
		/pgfplots/table/column name={\multicolumn{1}{c|}{##1}}}
},
display columns/14/.style={
	column name=ub,
},
display columns/15/.style={
	preproc/expr={100*##1},
	column name=gap,
		assign column name/.style={
		/pgfplots/table/column name={\multicolumn{1}{c|}{##1}}}
},
display columns/16/.style={,
	column name=$lb_{SDP}$,
},
display columns/18/.style={,
	column name=time,
},
display columns/17/.style={
	preproc/expr={100*##1},
	column name=gap,
},
display columns/19/.style={int detect,
	empty cells with={1},
	column name=rank,
},
every head row/.style={
	fixed,precision=2,
	before row={
		\caption{Detailed comparisons of experiments on instances with $n=400$}\label{tab:summary_400}\\
		\toprule  
		&  &  \multicolumn{6}{c|}{$\aleph = 5$}& \multicolumn{6}{c|}{$\aleph = 10$}& \multicolumn{6}{c}{$\aleph = 20$} \\
		& & \multicolumn{2}{c|}{Gurobi} &  \multicolumn{4}{c|}{Mosek} &\multicolumn{2}{c|}{Gurobi} & \multicolumn{4}{c|}{Mosek}& \multicolumn{2}{c|}{Gurobi} &  \multicolumn{4}{c}{Mosek}\\
	},
	after row={
		 & & & \% &  &  \% & & & &  \% & &  \% & & & &\% & & \% & &\\ 
		\endfirsthead
		\midrule} 
},
columns addvline/.style={
	every col no #1/.style={
	column type/.add={}{|}}
},
columns addvline/.list={1,3,7,9,13,15},
every row no 30/.style ={before row={\midrule}},
every row no 60/.style ={before row={\midrule}},
every last row/.style={after row=\bottomrule},
]{cardinality_summary_400.csv} 